\documentclass[12pt]{amsart}
\usepackage{amsmath, amssymb, latexsym}

\numberwithin{equation}{section}

\newcommand{\bd}{\mathbb{D}}
\newcommand{\bc}{\mathbb{C}}
\newcommand{\bt}{\mathbb{T}}

\renewcommand{\a}{\alpha}

\newcommand{\bnk}{{\beta_{nk}}}
\renewcommand{\l}{\lambda}
\newcommand{\lk}{{\lambda_k}}
\newcommand{\s}{\sigma}
\renewcommand{\r}{\rho}
\newcommand{\vro}{\varrho}

\renewcommand{\t}{\theta}
\renewcommand{\d}{\delta}
\newcommand{\dn}{{\delta_n}}
\newcommand{\dde}{\Delta}
\newcommand{\dd}{\Delta}
\renewcommand{\o}{\omega}
\newcommand{\on}{{\omega_n}}
\newcommand{\mn}{{\mu_n}}

\newcommand{\g}{\gamma}
\newcommand{\gn}{{\gamma_n}}
\newcommand{\gnk}{{\gamma_{nk}}}
\newcommand{\G}{\Gamma}
\newcommand{\ep}{\varepsilon}

\newcommand{\sh}{\#}
\newcommand{\di}{\operatorname{dist}}
\newcommand{\tr}{\operatorname{tr}\,}
\newcommand{\kker}{\mathrm{Ker\,}}
\newcommand{\rk}{\mathrm{rank\,}}
\newcommand{\nt}{\noindent}
\newcommand{\hf}{\frac{1}{2}}

\newcommand{\bsl}{\backslash}
\newcommand{\ovl}{\overline}
\newcommand{\prt}{\partial}
\newcommand{\ti}{\tilde}
\newcommand{\lb}{\left[}
\newcommand{\rb}{\right]}

\newcommand{\dsp}{\displaystyle}
\newcommand{\bv}{\bigvee}
\newcommand{\tl}{{\,\tilde{}}}

\newcommand{\cc}{\mathcal{C}}

\newcommand{\fs}{\mathfrak{S}}
\newcommand{\fd}{\mathfrak{D}}

\newcommand{\lrg}{\mathrm{(LRG)}}
\newcommand{\utb}{\mathrm{(UTB)}}
\newcommand{\utbm}{\mathrm{(UTB)}_*} 

\newcommand{\deo}{{\delta_0}}

\newcommand{\Ga}{\Gamma}
 
\newtheorem{corollary}{Corollary}[section]
\newtheorem{lemma}{Lemma}[section]

\newtheorem{theorem}{Theorem}[section]

\begin{document}

\title[Similarity to a normal operator]
{Operators similar to contractions and their similarity to a normal 
operator} 

\author{S. Kupin}

\date{November 5, 2001}

\address{Department of Mathematics, 253-37, Caltech, Pasadena CA 91125, USA.}

\address{Also: Math. Division, Institute of Low Temperature Physics, Lenin's ave., 
47, 61164 Kharkov, Ukraine.} 

\email{kupin@its.caltech.edu}

\maketitle

\vspace{-0.5cm}
\begin{abstract}
It has been established quite recently \cite{nb}, that a contraction $T$, 
having finite
defects ($\rk (I-T^*T),\; \rk(I-TT^*)<\infty$) and the spectrum $\s(T)$ not
filling in the closed unit disk $\ovl\bd$, is similar to a normal operator if
and only if
$$
\cc_1(T)=\sup_{\l\not\in\s(T)} ||(T-\l)^{-1}||\cdot\di(\l,\s(T))<\infty.\hspace{2cm} 
$$ 
It was shown later \cite{ntr}, that the ``only if" part of the assertion is no longer true if $T$ is not a contraction.

It is natural to ask what larger (than the class of contractions) classes of operators obey the linear resolvent growth test, formulated above. In this paper, we obtain results of the same type for operators that are similar to a contraction. Their proofs are based on a refined version of a theorem obtained  in \cite{ku1}. 
 
\end{abstract}  

\section*{Introduction}\label{s0}
It is common knowledge \cite{we} that an operator on a Hilbert space has an unconditionally convergent spectral decomposition if and only if it is similar to a normal operator. That is why it is important to have simple and efficient tests for similarity of a given operator to a normal one.

It has been established recently \cite{nb} that a contraction $T$ (i.e., an operator 
with norm less or equal to one), having finite defects ($\rk (I-T^*T),\;
\rk(I-TT^*)<\infty$) and  with the spectrum $\s(T)$ not filling in the closed
unit disk
$\ovl\bd$, is similar to a normal operator if and only if it has the so-called (LRG) property
$$
\cc_1(T)=\sup_{\l\not\in\s(T)} ||(T-\l)^{-1}||\cdot\di(\l,\s(T))<\infty.\eqno{\lrg}
$$
The abbreviation (LRG) stands for the linear resolvent growth. 
It was shown later \cite{ntr} that the ``only if" part of the assertion is no longer true if $T$ is not a contraction.

Therefore, it is natural to ask what larger (than the class of contractions) classes of operators obey the linear resolvent growth test for similarity to a normal operator, formulated above.  

Let us introduce some notation to formulate the main results of this paper. We say that a contraction $T$
has the uniform trace boundedness property (the (UTB) property, for short), if
$$
\cc_2(T)=\sup_{\mu\in\bd} \tr (I-b_\mu(T)^*b_\mu(T))<\infty,\eqno{\utb}
$$
where $b_\mu(T)=(I-\ovl\mu T)^{-1}(T-\mu)$ and $\mu\in\bd$.
Furthermore, an operator $L$ possesses the modified (UTB) property, if
$$  
\cc_2(L)=\sup_{\mu\in\bd} \tr |I-b_\mu(L)^*b_\mu(L)|<\infty.\eqno{\utbm}
$$
Here $|A|$ is the self-adjoint factor $(A^*A)^\hf$ of an operator $A$ in its polar decomposition.

The following theorem holds.
\begin{theorem} \label{th1}  
Let an operator $L$, acting on a Hilbert space, be similar to a contraction and   
$\s(L)\not=\ovl \bd$. If $L\in\lrg, L_c\in\utbm$, and $I-(L^*)_c^*(L^*)_c\in\fs_1$, the class of nuclear operators, then $L$ is similar to a normal operator.
\end{theorem}
Here $L_c$ stands for the complete part of $L$.
 
Theorem \ref{th1} is a product of a further development of an approach suggested in \cite{ku1}.
Nevertheless, the information provided by the latter work is far from being sufficient for the proof of the result. Essentially it will rely on the next theorem.

\begin{theorem} \label{th2}
Let $T$ be a contraction on a Hilbert space and 
$\s(T)\not =\ovl\bd$. If
$T\in\lrg$ and $T_c\in\utb$ along with $I-(T^*)_c^*(T^*)_c\in \fs_1$, 
then $T$ is similar to a normal operator.
\end{theorem}

We comment on the proofs of these theorems. One of the main ideas underlying the proof of Theorem \ref{th2} is that the (LRG) property for the original contraction $T$ allows one to split its part $T_1$ having spectrum on the unit circle $\bt$ from the complete part of $T$ with spectrum in $\bd$. These parts of $T$ are studied ``almost'' independently. For example, it turns out that the operator $T_1$ has to be similar to a unitary operator. Since the necessary and sufficient conditions for similarity to a unitary operator (see   Section~\ref{s12}) are valid for ge\-neral contractions, we do not impose any conditions on the defects of $T_1$. 

Nevertheless, since the above mentioned splitting uses the operator-valued corona theorem \cite{vas1}, its technical realization is quite complicated. 

Yet another important ingredient of the proofs is the following lemma on the conformally invariant nature of the (LRG) property.
\begin{lemma}[{\cite{ku1, ktr}}]\label{l01}
Let $T$ be a contraction on a Hilbert space. Then
$$
A_1\cc_1(T)\le\cc_1(b_\mu(T))\le A_2\cc_1(T),
$$
for all $\mu\in\bd$, $A_1$ and $A_2$ being absolute constants.
\end{lemma}

We turn to Theorem \ref{th1} now. Suppose that an operator $L$ satisfies its assumptions, and, in particular, it is similar to a contraction $T$. 
Surprisingly, the positivity condition $I-T^*T\ge 0$ ensures that
$\tr (I-T_c^*T_c)$ admits a bound in terms of $\tr |I-L_c^*L_c|$ in a ``regular" situation. 
The converse bound, generally speaking, is not true (see Lemma \ref{l3} and subsequent discussion). However, reducing the general situation to the ``regular" one requires some information on the sparseness of $\s(L)\cap\bd$. We obtain 
necessary facts combining the (UTB)$_*$ and the (LRG) properties. Again, Lemma \ref{l01} is implicitly present behind every step of the argument.

The paper is organized as follows. In Section \ref{s1}, we recall some well-known facts on the
Sz.-Nagy--Foia\c s functional model, geometry of point sets in $\bd$ and ideals of compact operators. Section \ref{s2} deals with the proof of Theorem \ref{th2}. Auxiliary facts needed for Theorem \ref{th1} are presented in Section \ref{s3}, and the theorem itself is proved in Section \ref{s4}. The section is concluded with corollaries of Theorem \ref{th1} and some open questions.

\smallskip
\nt{\it Acknowledgments.}
\ I am grateful to N. Nikolski and V. Vasyunin for helpful discussions on the problem.
I would like also to thank S. Albeverio and C. Foia\c s, who suggested Corollary \ref{c11} and Corollary  \ref{c3} as conjectures.

\section{Preliminaries}\label{s1}
The material presented in this section is of common knowledge and is
quoted only for the reader's convenience.

\subsection{}\label{s11}
We begin with recalling some standard notation. 

By default, all operators appearing in this paper live on separable Hilbert spaces.
The spectrum of an operator $A$ is denoted by $\s(A)$, and the point spectrum $\s_p(A)$ stays for $\{\l\in\s(A): \ \kker(A-\l I)\not=\{0\}\}$. An operator $A$ defined on
$H$ is called complete, if 
$$
H=\bv_{\l\in\s_p(A)}\kker (A-\l I)^{n(\l)},
$$
where $n(\l)$ is the  algebraic multiplicity of $\l\in\s_p(A)$. We put $H_c$ to be the latter linear span and we call $A_c=A|_{H_c}$ the complete part of $A$. The resolvent of $A$ is defined by the formula $R_\l(A)=(A-\l I)^{-1}$.

The operators $A_1\in L(H_1)$ and $A_2\in L(H_2)$ are called similar if there
exists a boundedly invertible operator $W\in L(H_2, H_1)$ such that
$A_2=W^{-1}A_1W$.
 
Let $E$ and $E_*$ be separable
Hilbert spaces. We denote by $L(E,E_*)$ the space of bounded linear
operators mapping $E$ into $E_*$. We put $L(E)=L(E,E)$. Furthermore, we
write $H^\infty(L(E,E_*))$ for the space of $L(E,E_*)$-valued
bounded analytic functions on the unit disk $\bd=\{z\in\bc: |z|<1\}$.
We put $BH^\infty(L(E,E_*))$ to be the unit ball of the space.
Similarly, we denote by $L^\infty(L(E,E_*))$ the space of the
$L(E,E_*)$-valued bounded measurable functions on the unit circle
$\bt=\{z\in\bc:  |z|=1\}$. As usual,  $L^2(E)=L^2(\bt, E)$ is
the Hilbert space of measurable on $\bt$ functions $f$, taking
values in $E$ such that 
$$
||f||^2=\frac{1}{2\pi}\int_0^{2\pi}
||f(e^{i\phi})||_E^2
d\phi<\infty,
$$ 
and $H^2(E)$ is the Hilbert space of $E$-valued
analytic in $\bd$ functions with the norm
$$
||f||^2=\sup_{0\le r<1}
\frac{1}{2\pi} \int_0^{2\pi}||f(re^{i\phi})||_E^2 d\phi<\infty.
$$

\subsection{}\label{s12}
A wide panorama of the subject we
discuss  here can be found in monographs \cite{nk, snf}.

We introduce some notation to define the functional model.
Let us fix a function $\t\in BH^\infty(L(E_*,E))$. We put
$\dde(t)=(I-\t(t)^*\t(t))^{1/2}\in
L^\infty(L(E_*))$, $0\le\dde(t)\le I \hbox{ a.e.~on }\bt$.

Further, we consider the so-called model space
\begin{equation}
K_\t=
\lb\begin{array}{c}
H^2(E)\\
\ovl{\dde L^2(E_*)}
\end{array}\rb
\ominus
\lb\begin{array}{c}
\t\\
\dde
\end{array}\rb H^2(E_*).\label{eq:model_1}
\end{equation}
We denote by $T_\t$ an operator defined on $K_\t$ by
the formula
$$
T_\t x=
\lb\begin{array}{c}P_+\ovl zx_1\\\ovl z x_2\end{array}\rb,
$$
where $P_+:L^2(E)\to H^2(E)$ is the Riesz orthogonal projection and 
$x=\lb\begin{array}{c} x_1\\x_2\end{array}\rb\in K_\t$. The operator is a contraction, $||T_\t||\le 1$,
and it is called the model operator.

Any contraction $T$, defined on a Hilbert space $H$, can be represented in the form $T=U_0\oplus T_0$,
where $U_0$ is a unitary operator and $T_0$ is a completely nonunitary
contraction ({\it c.n.u.~contraction}, to be brief), that is, none of the restrictions of the latter to its
reducing subspaces is unitary.

Moreover, for $\l\in\bd$, we define an operator-valued function $\t_T(\l)$ by the formula
\begin{equation}
\t_T(\l)=-T^*+\l D_{T}(I-\l T)^{-1}D_{T^*}|_{\fd_{T^*}},\label{e11}
\end{equation}
where $D_T=(I-T^*T)^{1/2},\ D_{T^*}=(I-TT^*)^{1/2}$,\; and \; $\fd_T=\ovl{D_T H},\ \fd_{T^*}=\ovl{D_{T^*}H}$.
The function $\t_T$ is called the characteristic function of $T$.
It can be shown that $\t_T\in BH^\infty(\fd_{T^*},\fd_T)$ and that
$\t_T$ is pure, that is, the only subspace $E\subset \fd_{T^*}$ where
$\t_T(t)|_E$ is a unitary constant a.e.~on $\bt$ is $E=\{0\}$.

The following theorem links the two series of definitions given above.
\newpage
\begin{theorem}[\cite{snf}]\label{thm:model}
\nt\begin{enumerate}
\item[\it i)] Any c.n.u.~contraction $T$ defined on a
Hilbert space $H$ is unitarily equivalent to $T_{\t_T}$, where
$\t_T\in BH^\infty(\fd_{T^*},\fd_T)$ is the characteristic function
of the contraction.
\item[\it ii)] Let $\t$ be a pure   function
from $BH^\infty(L(E_*,E))$. Then the contraction
$T_\t$ is completely nonunitary and its characteristic
function coincides with the
initial function $\t$.
\end{enumerate}
\end{theorem}
The characteristic function $\t_T$ of a contraction $T$ will usually be denoted by $\t$, and a 
c.n.u.~contraction itself will be identified with $T_{\t_T}$. 
It is transparent from formula \eqref{e11} that $\t_{T^*}=\t\tl$, where
$\t\tl(\l)=\t(\ovl\l)^*$. 

Finally, a contraction $T, \s(T)\subset\bt$, is similar to a unitary operator
if and only if the operators $\t(\l)$ are invertible for all $\l\in\bd$ and
$||\t(\l)^{-1}||\le C<\infty$ (see \cite{gk1}, \cite[ch.~9]{snf}).

\subsection{}\label{s14}
It is well-known that every invariant subspace $L$ of the model
operator
$T_\t$ ($T_\t L\subset L$) defines a certain regular factorization
$\t=\t_2\t_1$, \cite[ch.~7]{snf}. The converse is also true, that is,
every regular factorization of the characteristic function
$\t=\t_2\t_1$ of a contraction $T$ permits to
construct a $T$-invariant subspace $L_{\t_2}$. We set $T_{\t_2}=T|_{L_{\t_2}}$.

We refer to \cite[ch.~7]{snf} for the definition of regular factorizations, as
well as for their basic properties. For instance, it follows from the construction that the characteristic function of $\dsp
T_{\t_2}$ coincides with the pure part of $\t_2$. Also, it is not difficult to see that
factorizations $\t=\t_2\t_1$, where $\t_2$ is inner or
$\t_1$ is $*$-inner, are regular. 
In particular, when $\t_2$ (or $\t_1$) is two-sided inner, the model space $K_\t$ admits the following orthogonal decompositions
\begin{eqnarray}
K_\t&=&L_{\t_2}\oplus L_{\t_2}^\perp
=\left[\begin{array}{c} K_{\t_2}\\0\end{array}\right]
\oplus\left[\begin{array}{cc} \t_2& 0\\0&I\end{array}\right]K_{\t_1},\label{e32}\\
&&\nonumber\\
K_\t&=&L_{\t_2}\oplus L_{\t_2}^\perp=
{\lb\begin{array}{cc} I&0\\0& \t_1^*\end{array}\rb}K_{\t_2}\oplus
\left[\begin{array}{c}\t\\\dde\end{array}\right]{\t}^*_1K_{\t_1},\hspace{1cm}
\label{e31}
\end{eqnarray}
respectively. We recall that a function $\t\in H^\infty(L(E_*,E))$ is called inner
($*$-inner) if $\t(t)^*\t(t)=I \hbox{ a.e.~on }\bt$
($\t(t)\t(t)^*=I \hbox{ a.e.~on }\bt$). We say that $\t$ is two-sided inner if it is inner and $*$-inner. The function $\t$ is said to be outer ($*$-outer)
if\  $\ovl{\t H^2(E_*)}=H^2(E)$ (the function $\t\tl$ is
outer). 
 
\subsection{}\label{s15}
Let $T$ be a completely nonunitary
contraction defined on a Hilbert space $H$.
Let $\t=\t_2\t_1,\ \t_1\in BH^\infty(L(E_*,F))$, and $\t=\t'_2\t'_1,\ \t'_1\in BH^\infty(L(E_*,F'))$ be 
regular factorizations of $\t$, where $F, F'$ are some intermediate Hilbert spaces. Put $L=L_{\t_2},\ L'=L_{\t'_2}$ 
to be the corresponding invariant subspaces of $T$. Assume that $L+L'$ is dense in $H$.
We are interested in conditions for the angle between these subspaces
to be positive (or, in
other words, when $L\cap
L'=\{0\}$ and the sum $L+L'$ is closed).  

It is proved in \cite{vas1} that $L+L'$ is a direct decomposition of $H$
if and
only if the Bezout equation,
\begin{equation}
\Ga_1\t_1+\Ga'_1\t'_1=I,
\label{e1}
\end{equation}
is solvable with $\Ga_1\in H^\infty(L(F,E_*)),\ \Ga'_1\in
H^\infty(L(F',E_*))$, and an additional equation of the same type is solvable in
certain $L^\infty$ spaces as well. It is known (see references in
\cite{vas1}) that if $\t$ is two-sided inner, the sole
equation (\ref{e1}) is sufficient to have $H=L\dot{+}L'$. There
are some
other special cases, where the solvability of the equation
(\ref{e1}) implies the conclusion. The following
theorem, for instance,  is a result of the same type.
\begin{theorem}[\cite{nb},  Sect.~1.6]
\label{thm:cor_thm} Let $L,L'$ be invariant subspaces, defined by
regular factorizations $\t=\t_2\t_1,\t=\t'_2\t'_1$, and let the
sum $L+L'$ be dense in $H$. If $\t'_1$ is a $*$-inner function, the sum $H=L+L'$ is a direct sum
whenever equation \eqref{e1} is solvable.
\end{theorem}

\subsection{}\label{s16}
Recall that we may always factorize $\t$ as
$\t=\t_{inn}\t_{out}$ and  $\t =\t_{out*}\t_{inn*}$ \cite[ch.~5]{snf}, where
the function $\t_{inn}$ ($\t_{inn*}$) is inner
($*$-inner), and the function $\t_{out}$
($\t_{out*}$) is outer ($*$-outer), respectively. These are the so-called     
inner-outer and $*$-inner-outer factorizations.

We will use factorizations of the same type having more specific properties.
\begin{lemma}\label{l1}
Let $\t\in BH^\infty(L(E_*,E))$. Suppose that $\t=B\t_1=\t'_2B'$, where 
\begin{itemize}
\item[\it i)] the functions $B, B'$ are Blaschke-Potapov products {\rm (}see \cite{po}{\rm )}, admitting scalar multiples,
\item[\it ii)] the values of $\t_1$ and $\t'_2$ are invertible operators for every $\l\in\bd$ and, moreover, $||\t_1^{-1}||,||{\t'_2}^{-1}||\le C$ on $\bd$.
\end{itemize}
Then the sum 
$\dsp\left[\begin{array}{c} K_B\\0\end{array}\right]+
\lb\begin{array}{cc} I&0\\0&{B'}^*\end{array}\rb K_{\t'_2}$ 
is dense in $K_\t$.
\end{lemma}
The proof of the lemma and some comments on formula \eqref{e31} are quoted in the Appendix.

\subsection{}\label{s17}
Detailed information on the geometry of discrete sets in the unit disk can be
found in \cite{ga,nk}.

For $\l,\mu\in\bd$, we define
$\vro(\l,\mu)=|b_\mu(\l)|$ and 
$B_\d(\mu)=\{\l\in\bd: |b_\mu(\l)|\le\d\}$, where $0<\d<1$.
We say that the set
$\s=\{\l_k\}$ is sparse, if there
exists a $\d>0$ such that
\begin{equation}
B_\d(\l_1)\cap B_\d(\l_2)=\emptyset, \label{eq:de_carl_1}
\end{equation}
where $\l_1,\l_2\in\s$ and $\l_1\not=\l_2$.
The set $\s$ is called Carleson if
$$
\inf_{\mu\in\s} \prod_{\l\in\s\bsl\{\mu\}} |b_\mu(\l)|\ge\deo>0.
$$
We say that a set $\s$ is $N$-Carleson ($N$-sparse) if
it is a union of $N$ Carleson sets (of $N$ sparse
sequences), respectively.

\subsection{}\label{s18}
A nice reference on the subject of this subsection is \cite{gk1}.

Let $\fs_\infty$
denote the ideal of compact operators on $H$. The
Schatten--von-Neumann ideals, $\fs_p,\ 0<p\le\infty$, are defined
as
$$
\fs_p=\{A\in\fs_\infty:\sum_k s_k(A)^p<\infty\},
$$
where $s_k(A)=\l_k(A^*A)^{1/2}$ and $\l_k(A)$
are the eigenvalues of the operator $A$. 

If $A=A^*\in\fs_\infty$, we represent it as
$A=\sum_k\a_k(.,e_k)e_k$,
where $\{e_k\}_k$ are its normalized eigenvectors. We define $A_+$ and $A_-$ by the formulas
$$
A_+=\sum_{k:\; \a_k\ge 0}\a_k(.,e_k)e_k, \quad A_-=(-A)_+.
$$
Note that $A_+, A_-\ge 0$.
 
Furthermore, let $A\in\fs_1$ and $\{e_k\}_{k=1,\infty}$ be an arbitrary
orthonormal basis of $H$. It is known that the sum $\tr A=
\sum_k (Ae_k, e_k)$ converges and does not depend on the choice of
the orthonormal basis.

It is clear that if $A=A^*\ge 0$ and $A\in\fs_1$, then
$$
\tr A=\sum_j\l_j(A)=\sum_j s_j(A).
$$
This relation implies that
$\tr PAP\le \tr A$ for any orthogonal projection $P$ and
any
operator $A$ with the properties stated above. In particular, if
$k=\mathrm{rank\,}P<\infty$, then
\begin{equation}
k\min_{j=1,k} \l_j(PAP)\le \sum_{j=1}^k s_j(PAP)\le\tr PAP\le\tr A.
\label{eq:tr_class_1}
\end{equation}
 
\section{Contractions and Their Similarity to a Normal Operator}\label{s2}
 
\subsection{}\label{s22}
We begin this subsection by introducing a new technique. Let $\{\s_j\}_{j=1,N}$ 
be Carleson sets and $\s=\cup_j\s_j$.
Put 
$$
\vro_0=\min_{j=1,N} \ \inf_{z,w\in\s_j;\ z\not=w}\vro(z,w),
$$ 
and $\vro=\vro_0/4(N+1)$. Then consider a set given by $\cup_{z\in\s} B_\vro(z)$. The connectedness
components of the latter are denoted by $\{G_n\}$. The sets $\s^n$ are defined as $\s^n=\s\cap G_n$. The
construction of $\{\s^n\}$ implies that the sets have the following properties (see \cite[ch.~9]{nk}, \cite{vas2}, \cite{ku1}):
\begin{enumerate}
\item[\it i)] $\sh\s^n\le N$,
\item[\it ii)] $\max \{\vro(z,w):\ z,w\in\s^n\}\le \hf\vro_0$,
\item[\it iii)] $\vro(\s^n,\s\bsl\s^n)\ge\vro>0$,
\item[\it iv)] $\vro(\s^n,\prt G_n)=\vro$ \ and \ $\vro(\s,G_n)=\vro$.
\end{enumerate}
Suppose now that a function $\t\in BH^\infty(L(E_*,E))$ admits factorization $\t=\t'_2B'$, where
$B'$ is a Blaschke-Potapov product and operators $\t'_2(\l)\in L(F,E)$ are invertible for all $\l\in\bd$, $F$ being an intermediate space. 

\begin{lemma}\label{l2} 
Let $\t$ be as above and $\s={B'}^{-1}(\{0\})$ be an $N$-Carleson set.
Suppose that
$$
||\t(\l)^{-1}||\le C_1\sup_{\mu\in\s}\frac{1}{|b_\mu(\l)|}.
$$
Then $||{\t'_2}^{-1}||\le C_2$ on $\bd$.
\end{lemma}
 
\begin{proof}
Construct the above partition for the set $\s$.
It is obvious that the operator-valued function ${\t'_2}^{-1}=B'\t^{-1}$ is analytic on $G_n$ and $||{\t'_2}^{-1}||$ is subharmonic there. For $\l\in\prt G_n$, we have
$$
||{\t'_2}^{-1}(\l)||\le||\t^{-1}(\l)||\le\dsp C_1 \sup_{\mu\in\s}\frac{1}{|b_\mu(\l)|}\le \frac{C_1}{\vro},
$$
because of the property iv) of the partition of $\s$. By the maximum principle for subharmonic functions, we see that $||{\t'_2}^{-1}(\l)||\le C_1/\vro$ for $\l\in G_n$. 
The inequality also holds for $\l\in\bd\bsl\{\cup_n G_n\}$. Indeed, we have that
$\inf_{\mu\in\s}|b_\mu(\l)|\ge\vro$, and, as above,
$||{\t'_2}^{-1}(\l)||\le C_1/\vro$.
Hence, we see that $||{\t'_2}^{-1}||\le C_2$ on $\bd$, and the lemma is proved. 
\end{proof}

\subsection{}\label{s221}
We will use a result on the growth of subharmonic functions.
\begin{lemma}[{\cite[Sect.~23]{nkh}}]\label{l201} Let
$\s$ be a discrete subset of $\bd$ and let $u$ be a
subharmonic function on $\bc\bsl\s$ satisfying the inequality
$$
u(\l)\le \max\left\{\frac{1}{\di(\l,\s)},\frac{1}{|1-|\l||}\right\}.
$$
Then
$$
u(\l)\le\frac{A_0}{\di(\l,\s)},
$$
where $\l\in\bd\bsl\s$ and $A_0$ is an absolute constant.  
\end{lemma}

\subsection{}\label{s222}
In this subsection, we collect two more technical lemmas on contractions. The first one is proved in
\cite{nb}.

\begin{lemma}\label{l21}
Let $T\in\lrg$ be a contraction with $\s(T)\cap\bd$ being discrete.  Then the eigenvalues of \ $T$ are algebraically simple.
\end{lemma}

Below, $\t$ is the characteristic function of a c.n.u.~contraction $T$.
\begin{lemma}\label{l202}
\begin{itemize}
\item[\it i)] Operators $\t(\l)\in L(\fd_{T^*},\fd_T)$ are invertible for $\l\in\bd\bsl\s(T)$.
\item[\it ii)] If $T\in\lrg$ and $\l\in\bd\bsl\s(T)$, then  
$$
||\t^{-1}(\l)||\le C_3\sup_{\mu\in\s(T)\cap\bd}\frac1{|b_\mu(\l)|}.
$$ 
\end{itemize}
\end{lemma} 
The first claim of the lemma is proved in \cite[ch.~6]{snf}.
The second claim is \cite{ku1}, Lemma 3.2.
 
\subsection{}\label{s223}
A contraction $T$ is called a weak contraction if its defect operators $I-T^*T$ and $I-TT^*$ are of trace class and
$\s(T)\not=\ovl\bd$. Weak contractions have quite special properties. For instance,  their point spectrum is always a Blaschke sequence in $\bd$ and their  defect spaces coincide, that is, $U\fd_{T}=\fd_{T^*}$ for some unitary operator $U$.

\smallskip\nt
{\it Proof of Theorem \ref{th2}.}
Let $T$ be a contraction, acting on $H$, and satisfying the assumptions of the theorem.
Let $T=U_0\oplus T_0$ be its canonical decomposition (see Section~\ref{s12}). Obviously, $T$ is similar to a normal operator if and only if $T_0$ is. Therefore, we may assume from the beginning that the contraction $T$ is completely nonunitary. 

Let $\t\in BH^\infty(L(E_*,E))$ be the characteristic function of $T$.  
It can be readily seen that the spaces
$$
H_c=\bv_{\l\in\s(T)\cap\bd} \kker(T-\l I)^{n(\l)}\quad \mathrm{and}\quad
H'_c=\bv_{\l\in\s(T)\cap\bd} \kker(T^*-\ovl\l I)^{n'(\l)}
$$
are invariant with respect to $T$ and $T^*$ ($n(\l)$ and $n'(\l)$ stand for the algebraic multiplicities of 
root subspaces of $T$ and $T^*$). 
Consequently, the corresponding factorizations of $\t$ and $\t\tl$, say, $\t=B\t_1$ and $\t\tl=B'\tl\t'_2\tl$, are regular. We have that  $\t_{T_c}=B$ and $\t_{(T^*)_c}=B'\tl$. 

Furthermore, $I-T_c^*T_c, I-(T^*)_c^*(T^*)_c\in \fs_1$, and $\s=\s(T)\cap\bd\not=\bd$.  Hence, $T_c$ and $(T^*)_c$ are weak contractions, and their discrete spectra are Blaschke sequences. The functions  $B$
and $B'$ have scalar multiples and, clearly, are two-sided inner.
Moreover, the completeness criterion \cite{snf} says that they are Blaschke-Potapov products (see \cite{po}).

Lemma \ref{l21} yields that $T$ and $T^*$ do not have nontrivial root subspaces, that is, $n(\l)=n'(\l)=1$, for $\l\in\s$.
 
Now, put $\s_0=\s(T)\cap\bt$. Since $\di(\l,\s(T))=\min\{\di(\l,\s),$ $\di(\l,\s_0)\}$ and $\di(\l,\s_0)\ge |1-|\l||$, we get with the help of Lemma \ref{l201} that
$$
\begin{array}{rl}
||R_\l(T_c)||\le&\dsp||R_\l(T)||\le \cc_1(T)\max\left\{\frac1{\di(\l,\s)},\frac1{\di(\l,\s_0)}\right\}\\
&\\
\le&\dsp \frac{A_0\cc_1(T)}{\di(\l,\s)}.
\end{array}
$$
Consequently, $T_c\in\lrg$. This, together with $T_c\in\utb$, shows that $T_c$ is similar to a normal operator by \cite{ku1}, Theorem 1.1. A by-product of this conclusion is that $\s$ is an $N$-Carleson 
set for some integer $N$ (see \cite{ku1}, Corollary 3.3).

Furthermore, Lemma \ref{l202} implies that the operators $\t(\l)\in L(E_*,E)$ are invertible for $\l\in\bd\bsl\s$, and, moreover, we have the bound
$$
||\t^{-1}(\l)||\le C_3\sup_{\mu\in\s}\frac1{|b_\mu(\l)|}.
$$
Since the operators $B(\l)$ and $B'(\l)$ are invertible for $\l\in\bd\bsl\s$, the invertibility of $\t(\l)$ yields that the operators $\t_1(\l)$ and $\t'_2(\l)$ are invertible as well.  
Note that $||\t_1^{-1}||,||\t'_2{}^{-1}||\le C$ on $\bd$ as provided by Lemma \ref{l2}.
The last inequality implies, in particular, that $T_{\t'_2}=T|_{L_{\t'_2}}$ is similar to a unitary  operator (see Section \ref{s12}).

Lemma \ref{l1} claims that the sum
$$
\lb\begin{array}{c} K_B\\0\end{array}\rb+
\lb\begin{array}{cc} I&0\\0&{B'}^*\end{array}\rb K_{\t'_2}
$$ 
is dense in $K_\t$. One the other hand, we see that the equation
$$
\G'\t_2'+\G B=I
$$ 
is solvable with $\G'\in H^\infty(L(E_*,E))$ and $\G\in H^\infty(L(E))$. By Theorem \ref{thm:cor_thm}, the angle between  the subspaces $H_c$ and $L_{\t'_2}$ is strictly positive and the sum $H_c+L_{\t'_2}$ is closed. Consequently, an operator orthogonalizing these subspaces is bounded and has a bounded inverse.  

Hence, the operator $T$ is similar to $T|_{H_c}\oplus T|_{L_{\t'_2}}$, which is similar in turn to the orthogonal sum of a diagonal operator and a unitary operator. This means that $T$ is similar to a normal operator, and the theorem is proved.  \hfill $\Box$

\medskip
It might seem that our assumptions on $T$ are quite asymmetric. Indeed, we require that $T_c\in\utb$ and
$I-(T^*)_c^*(T^*)_c$ be of the trace class only. This point is explained by the following lemma, which we quote without proof.

\begin{lemma} Let $T\in\lrg$ be a completely nonunitary contraction and $\s(T)\not=\ovl\bd$. Then the following assertions are equivalent:
\begin{itemize}
\item[\it i)] $T_c, (T^*)_c\in\utb$,
\item[\it ii)] $I-(T^*)_c^*(T^*)_c\in\fs_1$ and $T_c\in\utb$,
\item[\it iii)] $I-T_c^*T_c\in\fs_1$ and $(T^*)_c\in\utb$,
\item[\it iv)] $T_c\in\utb$ {\rm(}or $(T^*)_c\in\utb${\rm)} and $\s(T)\cap\bd$ is a discrete set.
\end{itemize}
\end{lemma}
 
\section{Some Auxiliary Propositions}\label{s3}

\subsection{}\label{s30}
The proof of Theorem \ref{th1} relies on several auxiliary propositions. 
Their proofs are close in spirit to those of \cite{ku1} (see Lemma 3.4 and below).  

Let $L$, defined on $H$, be an operator similar to a c.n.u.~contraction $T$ on $H_0$.   
In this section, we suppose that $T,L\in\lrg$, and therefore the operators do not have nontrivial root subspaces.
It is clear that $\s_p(L)=\s_p(T)$ and, for instance, $\s(L)\cap\bd=\s_p(L)\cup\s_p(L^*)$. We put
$\s_p(L)=\{\l_k\}$ and denote by $\{X_\l\}_{\l\in\s_p(T)}, \
X_\l=\kker (T-\l I)$, and  
$\{Y_\l\}_{\l\in\s_p(L)}, \ Y_\l=\kker (L-\l I)$, the families of the eigenspaces 
of $T$ and $L$, respectively. The subspaces $H_{0c}$ and $H_c$ are defined by relations
$$
H_{0c}=\bv_{\l\in\s_p(T)}\kker(T-\l I),\quad H_{c}=\bv_{\l\in\s_p(L)}\kker(L-\l I).
$$

\begin{lemma}\label{l30}
Let $L$ be an operator similar to a c.n.u.~contraction and $L_c\in\utbm$. Then
there exists an integer $M_1$ such that $\dim Y_\l\le M_1$
for every $\l\in\s_p(L)$.
\end{lemma}
\begin{proof}
Define an integer $k$ as $\dim Y_\l=\dim\kker b_\l(L), \ \l\in\s_p(L)$. Denote by $P_\l: H\to Y_\l$ the
orthogonal projection to $Y_\l$. Observing that $(I-b_\l(L)^*b_\l(L))|_{Y_\l}=I$, we get with the help of  \eqref{eq:tr_class_1}
$$
k=k\,\min_{j=1,k} \l_j(P_\l|I-b_\l(L_c)^*b_\l(L_c)|P_\l)\le \tr |I-b_\l(L_c)^*b_\l(L_c)|\le \cc_2(L_c).
$$
The lemma is proved.
\end{proof}

\subsection{}\label{s31}
We will prove in this subsection that if $L$ satisfies the assumptions of Theorem \ref{th1}, its point spectrum is not ``thicker" than a Blaschke sequence. 
\begin{lemma}\label{l3}
Let an operator $L\in\lrg$ be similar to a contraction $T$, $I-L_c^*L_c\in\fs_1$, and $\s(L)\not=\ovl\bd$. If   $\s_p(L)\cap B_\d(0)=\emptyset$ for some $\d>0$. Then
$$
\tr (I-T_c^*T_c)\le C_4(\d)\tr |I-L_c^*L_c|.
$$ 
\end{lemma}
\begin{proof}
Let $V: H\to H_0$ be the operator intertwining $L$ and $T$. For any integer $n$, we put
$$
H_{0n}=\bv_{k=1,n} X_\lk, \quad H_n=\bv_{k=1,n} Y_\lk.
$$
It is clear that $VH_n=H_{0n}$, and $N_n=\dim H_n=\dim H_{0n}< \infty$ by Lemma \ref{l30}. We define,
further,
$$
T_n=T|_{H_{0n}}, \quad L_n=L|_{H_n}
$$
and $V_n=V|_{H_n}: H_n\to H_{0n}$. It follows immediately from the definitions that $L_n=V_n^{-1}T_nV_n$. We denote by $\bnk$ and $\gnk$ the eigenvalues of operators $L_n^*L_n$ and $T_n^*T_n$, correspondingly.
We note that
$$
\begin{array}{rl}
\prod_{k=1}^{N_n}\bnk=&\det L_n^*L_n=|\det L_n|^2=|\det V^{-1}_nT_nV_n|^2\\
=&\det T_n^*T_n=\prod_{k=1}^{N_n}\gnk.
\end{array}
$$
Furthermore, taking $\l=0$ in the (LRG) inequality for $L$, we obtain
$$
||(L_n^*L_n)^{-1}||\le\frac{\cc_1(L)^2}{\di(0,\s(L))^2}\le\left(\frac{\cc_1(L)}{\d}\right)^2.
$$
Setting $C_4(\d)=(\cc_1(L)/\d)^2$, we get $1/\bnk\le C_4(\d)$.
  
Applying inequality $\log 1/x\ge 1-x, x>0$, we see that
$$
\log\frac1{\prod_{k=1}^{N_n} \gnk}\ge\sum_{k=1}^{N_n}(1-\gnk).
$$
On the other hand,
$$
\begin{array}{rl}
\dsp\log \frac1{\prod_{k=1}^{N_n}\gnk}=&\dsp\log\frac1{\prod_{k=1}^{N_n}\bnk}\le \sum^{N_n}_{k:\ \bnk\le 1}\log\frac1{\bnk}\\
&\\
\le&C_4(\d)\sum_{k:\ \bnk\le 1}^{N_n}(1-\bnk)\le C_4(\d)\sum_{k=1}^{N_n} |1-\bnk|,
\end{array}
$$
where we have used that $\log 1/x\le 1/x-1$ for $x>0$. Hence, we have
$$
\tr (I-T_n^*T_n)\le C_4(\d)\tr |I-L_n^*L_n|\le C_4(\d)\tr |I-L_c^*L_c|.
$$
Passing to the limit by $n\to\infty$, we obtain the conclusion of the lemma.
\end{proof}

It turns out that there is no converse bound to that one obtained in Lemma \ref{l3}. Indeed,  an analysis of the proof shows that 
$\tr (I-T_c^*T_c)$ controls only the difference
$|\tr (I-L_c^*L_c)_+ - \tr (I-L_c^*L_c)_-|$, whence 
$$
\tr |I-L_c^*L_c|=\tr (I-L_c^*L_c)_++\tr (I-L_c^*L_c)_-,
$$ 
see  Section~\ref{s18} for notation.

The inclusion $L\in\lrg$ entails $T\in\lrg$, and, by virtue of Lemma \ref{l01}, $b_\mu(T)\in\lrg$.
Consequently, $b_\mu(L)\in\lrg$, and moreover, 
$$
A'_1\cc_1(L)\le\cc_1(b_\mu(L))\le A'_2\cc_1(L)
$$
with constants $A'_1$ and $A'_2$ not depending on $\mu\in\bd$.
\begin{corollary}\label{c1}
If $L\in\lrg$ is similar to a c.n.u.~contraction, $L_c\in\utbm$, and $\s(L)\not=\ovl\bd$, then $\s_p(L)$ is a Blaschke sequence.
\end{corollary}
\begin{proof}
Pick a point $\mu\in\bd\bsl\s(L)$ and $\d>0$ such that $\s(L)\cap B_\d(\mu)=\emptyset$. Consider an operator given by
$b_\mu(L)=V^{-1}b_\mu(T)V$. By Lemma \ref{l01}, $L\in\lrg$ implies that $b_\mu(L)\in\lrg$. Lemma \ref{l3}, applied to $b_\mu(T)$ and $b_\mu(L)$, shows that
$$
\tr (I-b_\mu(T_c)^*b_\mu(T_c))\le C_4(\d)\tr |I-b_\mu(L_c)^*b_\mu(L_c)|\le C_4(\d)\cc_2(L_c),
$$
or $I-b_\mu(T_c)^*b_\mu(T_c)\in\fs_1$ and this suggests that $b_\mu(T_c)$ is a weak contraction. Hence, 
$\s_p(T_c)$ is a Blaschke sequence.
\end{proof} 

\subsection{}\label{s32}
Now, we may sharpen the conclusion of the previous subsection.

\begin{lemma}\label{l4}
Let $L$ be an operator as in Corollary \ref{c1}. Then
there exists a $\d_0>0$ such that 
$\sh B_\d(\mu)\cap\s_p(L)\le M_2$ for every $0<\d<\d_0$ and every $\mu\in\bd$.
\end{lemma}
\begin{proof}
Suppose that the claim of the lemma is false and there exist sequences $\{N_n\}$ and $\{\d_n\}, \d_n>0$, such that $N_n\to\infty$ and $\d_n\to 0$, when $n\to \infty$. Suppose also that we may find a sequence $\{\mu_n\}\subset \bd$ with the property $\sh B_{\d_n}(\mu_n)\cap\s_p(L)\ge N_n$. Let $\o_n$ be a subset of $B_{\d_n}(\mu_n)\cap\s_p(L)$. Denote $\sh\o_n$ by $N'_n$; the choice of $\o_n$ will be made precise later.

Now, consider the subspace $Y_\on=\bigvee_{\l\in\on} Y_\l$ and the operator $L_\on=L|_{Y_\on}$. It is not difficult to see that $b_\mu(L_\on)=b_\mu(L)|_{Y_\on}$ and $\s_p(b_\mu(L_\on))=b_\mu(\on)$ for any $\mu\in\bd$.

We estimate the norm of $b_\mn(L_\on)$ with the help of the Riesz-Dunford calculus formula \cite{d1}. To do that, we put $\g=\{z:|z-\l|=\ep, \l\in b_\mn(\on)\}$ and we choose $0<\ep\le\d/2$ sufficiently small to guarantee that different circles composing $\g$ do not intersect. We have the following formula for the operator $b_\mn(L_\on)$  
$$
b_\mn(L_\on)=-\frac{1}{2\pi i}\int_{\g_n}zR_z(b_\mn(L_\on))\,dz.
$$
Its norm may be estimated as follows
$$
\begin{array}{rl}
||b_\mn(L_\on)||\le&\dsp\frac{1}{2\pi}\int_\gn |z|\cdot ||R_z(b_\mn(L_\on))||\,|dz|\\
&\\
\le&\dsp\frac{|\gn|}{2\pi}\sup_{z\in\gn}|z|\sup_{z\in\gn}||R_z(b_\mn(L_\on))||
\end{array}
$$
$$
\begin{array}{rl}
\hspace{2.5cm}\le&\dsp C_5\, \ep\, \sh\on\, \d_n\cdot\frac{\cc_1(L)}{\ep}=C_6 N'_n\d_n,
\end{array}
$$
where $C_6$ does not depend on $n$.
We have used here that, if $L\in\lrg$, then $b_\mu(L)\in\lrg$ by observation quoted before Corollary \ref{c1}.

Since $N_n\to \infty$ and $\d_n\to 0$, we may choose $\on$ by dropping a right number of points from $B_\dn(\mn)\cap\s_p(L)$ in such a way that $N'_n=\sh \on\le 1/(2C_6\d_n)$, and $N'_n\to\infty$. This choice of $\on$ immediately gives that $||b_\mn(L_\on)||\le\hf$.

Now, estimate $\tr |I-b_\mn(L_\on)^*b_\mn(L_\on)|$ from below with the help of inequality \eqref{eq:tr_class_1}
$$
\frac{3}{4}N'_n\le \tr (I-b_\mn(L_\on)^*b_\mn(L_\on))=\tr |I-b_\mn(L_\on)^*b_\mn(L_\on)|.
$$
On the other hand,
$$
\tr |I-b_\mn(L_\on)^*b_\mn(L_\on)|\le \tr |I-b_\mn(L_c)^*b_\mn(L_c)|\le \cc_2(L_c), 
$$
and since $N'_n\to\infty$, we get a contradiction. The lemma is proved.
\end{proof}
 
\section{Proof of Theorem \ref{th1} and Its Corollaries}\label{s4}
\subsection{}\label{s41}

{\it Proof of Theorem \ref{th1}.}\quad 
Let an operator $L$, acting on $H$, satisfy the assumptions of the theorem. Assume that $L$ is similar to a contraction, defined on $H_0$.
Without loss of generality, we suppose that  $T$ is completely nonunitary. Indeed, if $L$ is similar to a normal operator, then so is $T$. On the other hand, $T$ can be represented as $T=T_0\oplus U_0$, where $T_0$ is a c.n.u.~contraction and $U_0$ is a unitary operator (see Section \ref{s12}). We see that $T$ is similar to a normal operator if and only if $T_0$ is. Hence, we may assume that $T=T_0$.

Since $\s(L)\not=\ovl\bd$, Corollary \ref{c1} shows that $\s_p(T)$ satisfies the Blaschke condition.
Since $T\in\lrg$, the eigenvalues of $T$ are algebraically simple by Lemma \ref{l21}.

We prove that $T$ is similar to a normal operator with the help of Theorem \ref{th2}. 

To apply it, we prove first that $I-(T^*)_c^*(T^*)_c\in\fs_1$.
Take $\mu\in\bd\bsl\s(T)$ and consider $L'=b_\mu(L)^*|_{H'_c}$, where $H'_c=\bv_{\l\in\s_p(L^*)}\kker(L^*-\l I)$. A simple computation \cite[ch.~6]{snf} gives that
\begin{equation}
I-{L'}^*L'=S^*(I-(L^*)_c^*(L^*)_c)S,\label{e3}
\end{equation}
where $S$ is a bounded invertible operator. Hence, if  $I-(L^*)_c^*(L^*)_c\in\fs_1$, then $I-{L'}^*L'\in\fs_1$, and, in light of Lemma \ref{l3}, 
$I-{T'}^*T'\in\fs_1$, where $T'=b_{\ovl\mu}((T^*)_c)$. 
Using \eqref{e3} with $T'$ and $(T^*)_c$ instead of $L'$ and $(L^*)_c$, we see that $I-(T^*)_c^*(T^*)_c\in\fs_1$.

It remains to show that
\begin{equation}\label{eq1}
\cc_2(T_c)=\sup_{\mu\in\bd}\tr (I-b_\mu(T_c)^*b_\mu(T_c))<\infty.
\end{equation}
This inequality is proved in the lemma below. \hfill$\Box$
\begin{lemma}\label{l7}
If $L\in\lrg$ is an operator similar to a contraction $T$, $L_c\in\utbm$, 
and $\s(L)\not=\ovl\bd$, then $T_{c}\in\utb$.   
\end{lemma}
\begin{proof}
It was already mentioned that $\s_{p}(T)=\s_{p}(L)$ is a Blaschke sequence 
in $\bd$. Moreover,     
Lemma \ref{l30}, applied to the operator $L$, gives that $\dim Y_\l\le 
M_1, \ \l\in\s_{p}(L)$, and $\s_p(L)$ is $N$-sparse with certain $N$ by Lemma \ref{l4}.
Now, taking
$V_{c}=V|_{H_{c}}: H_c\to H_{0c}$ (see Section~\ref{s30} for notation), we  get that $L_c=V_c^{-1}T_cV_c$.
We take $\d>0$, given by Lemma \ref{l4}, and define
$$
\s_1=b_\mu(\s_p(T))\cap B_\d(0),\quad \s_2=b_\mu(\s_p(T))\bsl\s_1.
$$ 
Set also $H_{0j}=\bv_{\l\in\s_{j}}X_\l$ and consider operators 
$b_\mu(T)\vert_{H_{01}}$ and $b_\mu(T)\vert_{H_{02}}$. 
Observe that the operator $b_\mu(T)\vert_{H_{02}}$ has the (LRG) property by \cite{ku1}, Lemma 3.7.
Put $N_{1}=\dim H_{01}$.
Then, choose an orthonormal basis $\{e_k\}_{k=N_1+1,\infty}$ in the subspace 
$H_{02}$,   and complete it up to a basis of the whole space $H_{0c}$. 
We have, by Lemma \ref{l3},
$$
\begin{array}{rl}
\tr (I-b_{\mu}(T_c)^{*}b_{\mu}(T_{c}))=&
\sum_{k=1}^\infty ((I-b_{\mu}(T_c)^{*}b_{\mu}(T_{c}))e_k,e_k)\\
&\\
\le&N_1+\sum_{k=N_1+1}^\infty ((I-b_{\mu}(T_c)^{*}b_{\mu}(T_{c}))e_k,e_k)\\
&\\
\le&N_1+\tr P_{H_{02}}(I-b_{\mu}(T_c)^{*}b_{\mu}(T_{c}))P_{H_{02}}\\
&\\
\le& M_1M_2+C_4(\d)\cc_{2}(L),
\end{array}
$$ 
where $P_{H_{02}}: H_{0c}\to H_{02}$ is the orthogonal projection from 
$H_{0c}$ on $H_{02}$ and $C_4(\d)$ is the constant from Lemma 
\ref{l3}, computed for $b_\mu(T)|_{H_{02}}$.
Inequality \eqref{eq1} is proved.
\end{proof}  

\subsection{}
We  mention here some corollaries of Theorem \ref{th1}.
\begin{corollary}\label{c11}
Let $L\in\lrg\cap\utbm$ be an operator similar to a contraction and  
$\s(L)\not=\ovl\bd$. Then it is similar to a normal operator.
\end{corollary}
The proof is immediate from the inequality
$$
\tr |I-b_\mu(L_c)^*b_\mu(L_c)|\le \tr |I-b_\mu(L)^*b_\mu(L)|
$$
and a  remark that $\tr |I-(L^*)_c^*(L^*)_c|\le\tr|I-LL^*|<\infty$.
\begin{corollary}\label{c2}
Let $L$ be an operator with finite defects and similar to a contraction. 
Suppose that  
$\s(L)\not=\ovl\bd$.
It is similar to a normal operator if and only if $L\in\lrg$.
\end{corollary}
The corollary is self-evident, since $\cc_{2}(L_c),\cc_2((L^*)_c)\le \rk \vert I-L^{*}L\vert$.
 
Recall that a bounded operator $L$ on a Hilbert space $H$ has a $\r$-dilation, $\r>0$, if there 
exists a Hilbert space $\ti H, \ H\subset \ti H$ and a unitary operator 
$U$ on it with the property
$$
L^n=\r P_H U^n,
$$
where $P_H:\ti H\to H$ is the orthogonal projection. It is well-known \cite[ch.~2]{snf} 
that if $L$ has a $\r$-dilation, then it is similar to a contraction. This fact yields the following corollary.
\begin{corollary}\label{c3}
Let $L$ be an operator having a $\r$-dilation and $\s(L)\not =\ovl\bd$. 
If $L\in\lrg$, $L_c\in\utbm$, and $I-(L^*)_c^*(L^*)_c\in\fs_1$, then
$L$ is similar to a normal operator.
\end{corollary}

It is reasonable to put questions on similarity 
to a normal operator for wider 
classes of operators. For instance, we do not know whether the (LRG)-type criteria 
work for polynomially bounded operators (power bounded operators, operators with spectral radius less or equal to one). Neither do we know what techniques should be applied in studying 
the question in this new setting.

\section*{Appendix}
Let $\t=B\t_1=\t'_2B'$ be regular factorizations of $\t$ described in Lemma \ref{l1}. It follows directly from its assumptions that $B$ and $B'$ are two-sided inner functions. 

Consider the second factorization in more detail. Define a map 
$Z':\ovl{\dd L^2(E_*)}$ $\to$ $\ovl{\dd'_2 L^2(E_*)}$ by the formula $Z'(\dd g)=\dd'_2B'g$, where $g\in L^2(E_*)$ and $\dd'_2=(I-{\t_2'}^*{\t_2'})^{1/2}$.  Since the factorization $\t=\t'_2B'$ is regular, the map $Z'$ is unitary \cite{snf}.
Moreover, we have that $\dd={B'}^*{\dd'}_2B'$ on $\bt$, and ${Z'}^{-1}(\dd'_2g)=\dd {B'}^*g={B'}^*\dd'_2g$.   Consequently, we obtain (see \cite{snf}, Theorem 1.1, ch.~7)
$$
L_{\t'_2}=\lb\begin{array}{cc}I&0\\0&{Z'}^{-1}\end{array}\rb K_{\t'_2}=
\lb\begin{array}{cc}I&0\\0&{B'}^*\end{array}\rb K_{\t'_2},
$$
and
$$
L^\perp_{\t'_2}=\lb\begin{array}{c}\t'_2\\{Z'}^{-1}{\dd'}_2\end{array}\rb H^2(E_*)
\ominus\lb\begin{array}{c}\t\\\dd\end{array}\rb H^2(E_*)=
\lb\begin{array}{c}\t \\\dd\end{array}\rb {B'}^*K_{B'},
$$

\medskip\nt
(see Sections~\ref{s14}, \ref{s16} for notation). This is exactly formula \eqref{e31}.

\medskip
\nt{\it Proof of Lemma \ref{l1}.} 
Suppose that $f\in K_\t$ is orthogonal to the sum appearing in Lemma \ref{l1}. We want to prove that, indeed, $f=0$.

It follows from the assumptions of the lemma and orthogonal decompositions \eqref{e32}, \eqref{e31}, that
$$
f=\left[\begin{array}{c}f_1\\f_2\end{array}\right]
\in \left[\begin{array}{cc} B&0\\0&I\end{array}\right]K_{\t_1}\cap
\left[\begin{array}{c}\t\\ \dde\end{array}\right]{B'}^*K_{B'}.
$$
Hence, we have $Bg_1=\t'_2g_2$ for some $g_1\in H^2(E)$ and $g_2\in K_{B'}$. Using the equality $\t=B\t_1=\t'_2B'$, we get $g_1=\t_1{B'}^{-1}g_2$. Since functions $B$ and $B'$ have scalar multiples, their inverses $B^{-1}, {B'}^{-1}$ are meromorphic on $\bd$ and their only singularities are poles forming
a Blaschke sequence. Furthermore, since $||\t_1^{-1}||\le C$ on $\bd$, the function 
$\t_1^{-1}g_1$  lies in $H^2(E_*)$, and so does ${B'}^{-1}g_2$.  
Consequently, $g_2=B'g_3$ with $g_3\in H^2(E_*)$. On the other hand, $g_2\in K_{B'}$, and we have 
$g_2=g_1=0$, and hence $f=0$. \hfill$\Box$

\end{document}